\documentclass[12pt]{amsart}
\usepackage{amssymb}
\newtheorem{theorem}{Theorem}

\newtheorem{coro}[theorem]{Corollary}
\newtheorem{lemma}[theorem]{Lemma}

\newtheorem{assum}{Assumption}

\newcommand\mcS{\mathcal{S}}

\newcommand\D{\mathcal{D}}
\newcommand\C{\mathbf{C}}
\newcommand\N{\mathbf{N}}
\newcommand\R{\mathbf{R}}

\newcommand\dd{\,{\rm d}}
\newcommand{\cdo}{\, \cdot \,}
\newcommand\wrt{\,{\rm d}}
\newcommand\ld{L^2(TX)}
\newcommand\lj{L^1(TX)}

\newcommand\gr{\nabla}

\def\d{\rho}
\def\p{{\partial}}
\def\cprime{$'$}
\def\re{{\mbox{\tiny\rm  Re}\,\,}}

\newcommand\mua{\mu(B(x,r))}

\begin{document}

\setcounter{section}{0}
\renewcommand{\theenumi}{\alph{enumi}}
\renewcommand{\labelenumi}{\textrm{(\theenumi)}}
\numberwithin{equation}{section}

\title[Riesz transform, Gaussian bounds and the method of wave equation]
{Riesz transform, Gaussian bounds and the method of wave equation}

\author{Adam Sikora}
\address{Adam  Sikora, Department of Mathematical Sciences, 
New Mexico State University, Las Cruces, NM 88003-8001, USA}
\email{asikora@nmsu.edu}

\def\d{\rho}
\def\C{{\bf C}}
\def\p{{\partial}}
\def\N{{\bf N}}
\def\R{{\bf R}}
\def\e{\varepsilon}
\def\q{{q'}}
\newcommand{\supp}{\operatorname{supp}}

\begin{abstract}
For an abstract self-adjoint operator $L$ and a local operator $A$ we study
the boundedness of the Riesz transform $AL^{-\alpha}$ 
on $L^p$  for some $\alpha >0$.
A very simple  proof of the obtained result is  based on the 
finite speed propagation property for the solution of the
corresponding wave equation. We also discuss the relation between the
Gaussian bounds and the finite speed propagation property. 
Using the wave equation methods we obtain
a new natural form of the Gaussian bounds for the heat kernels for 
a large class of the generating operators. 
We describe a surprisingly elementary proof of the finite speed propagation 
property in a more general setting than it is usually considered in the
literature.

As an application of the obtained results  we prove boundedness of the 
Riesz transform on  $L^p$ for all $p\in (1,2]$ for Schr\"odinger
operators with positive potentials and electromagnetic fields. 
In another application we
discuss the Gaussian bounds for the Hodge 
Laplacian and  boundedness of the Riesz transform on $L^p$ of the Laplace-Beltrami
operator on Riemannian manifolds  for $p>2$ . 

\end{abstract}

\maketitle

\section{Introduction}
Let $\Delta=-\sum_{i=1}^{n}\partial_i^2$ be the standard Laplace operator
 acting on
$\R^n$. Then the corresponding Riesz transform is defined by the formula  
$\partial_j \Delta^{-1/2}.$
The $L^p$ continuity of the   Riesz transform for all $p\in (1,\infty)$
is one of the most 
important and celebrated results in
analysis.
Papers devoted to the study of the Riesz transform 
and its generalizations are too
numerous to list here. Hence we would like to mention only
a few most relevant works 
\cite{Ale, Ba, CD1, CD2, ERS2, ERS1, GS, Li, Lo, Sal, Sh1, Sh2, Sj, St1, Str}.

 The operator $\gr L^{-1/2}$,
where $\gr$ is the gradient  and $L$ is the Laplace-Beltrami operator 
on a  Riemannian
manifold $M$, is a natural generalization of the classical
Riesz transform.  $L^2$ boundedness of 
the Riesz transform $\gr L^{-1/2}$ is a consequence of
the equality $\|\gr f\|_{L^2}=\|L^{1/2}f\|_{L^2}$, which is actually 
the definition of the operator $L$. 
In \cite{Str} Strichartz asked whether one could extend $L^p$ continuity
of the classical Riesz transform to the setting of 
Laplace-Beltrami operators described above. 
An answer to this question was given  
in \cite{CD1} for $1\le p \le 2$.
In \cite{CD1} Coulhon and Duong proved that if
$M$ is a complete Riemannian manifold which satisfies the doubling volume
property (see Assumption~\ref{1a}), $L$ is the Laplace-Beltrami 
operator on $M$ 
and the heat kernel corresponding to $L$ satisfies the Gaussian bounds
then the Riesz transform $\gr L^{-1/2}$ is of weak type $(1,1)$
and so bounded on $L^p$ for
all $p\in (1,2]$. Note that 
$(\partial_j\Delta^{-1/2})^{*}=-\partial_j\Delta^{-1/2}$
so the boundedness of the standard Riesz transform
$\partial_j\Delta^{-1/2}$ for $p\in (1,2]$ implies continuity 
of the standard Riesz transform for all $p\in (1,\infty)$. 
Surprisingly in general  the Riesz transform corresponding to the Laplace-Beltrami operator $\gr L^{-1/2}$ is no
longer necessarily continuous for $p>2$ even under the above assumptions 
(see \cite{CD1} for a counterexample).

One of the main aims of this paper is 
to remove any assumptions about the nature
of the operator $L$ from the result obtained in \cite{CD1}. In 
Theorem~\ref{gl} below we 
consider an abstract self-adjoint positive definite operator.
We show  that if an operator $L$ satisfies finite speed propagation property 
for the solutions of the corresponding wave equation,
  $A$ is a local
operator    and $AL^{-\alpha}$ is bounded on $L^2$ for some $\alpha >0$ then
the operator $AL^{-\alpha}$ is automatically bounded
on all $L^p$ for $p\in (1,2]$ and  of weak type $(1,1)$. Thus, it turns out
that one does not have to assume that $L$ is the Laplace-Beltralmi operator in
\cite{CD1} and 
that the finite speed propagation property is the only essential 
assumption in \cite{CD1}.
Removing assumptions about the nature of the operator $L$ allows us to
 study the Riesz transform for Schr\"odinger
operators with positive potentials and electromagnetic fields. 
Such Riesz transforms were investigated in \cite{Li, CDE, Sh1, Sh2}.
Results which we describe here generalize and strengthen  a part 
of the results described in  \cite{Li, CDE, Sh1, Sh2}.

We start our discussion with a description of the equivalence of
the finite speed propagation property and the $L^2$ version of the 
Gaussian bounds (see condition~(\ref{fs2})). This allows us to obtain a very
elegant and straightforward proof of the  finite speed propagation property  
(see Theorem~\ref{fsp}) in a
 more general setting than usually seen in the literature (see for example 
\cite{Fo, Ta}). Theorem~\ref{fsp}
describes a large and natural class of examples of possible applications
of our main results. 
 
Another goal of this paper is to prove that on-diagonal
bounds for the heat kernel and finite speed propagation property imply
off-diagonal Gaussian bounds (see Theorem~\ref{gl2}). 
It is well known that for the Laplace-Beltrami operator, for 
the sublaplacians acting on Lie groups, and more generally for diffusion
semigroups on-diagonal bounds imply Gaussian off-diagonal bounds
(see for example \cite{Co, DP, Gr1}).
The advantage of our approach is that Theorem~\ref{gl2} again
holds without any assumptions about the nature of the semigroup
generator. Our only assumption is the finite propagation speed. 
And so for example Theorem~\ref{gl2} together with Theorem~\ref{fsp} 
show that on-diagonal estimates 
imply sharp Gaussian off-diagonal estimates for the
heat kernels generated by the Hodge Laplacian acting on $p$-forms
(see Corollary~\ref{co}).
The proofs from \cite{Co, DP, Gr1}
do not easily  generalize  to the 
De Rham-Hodge Laplacian setting. It turns out that the understanding of  
the behavior of 
the heat kernel generated by the Hodge Laplacian is a useful tool in the
study of the $L^p$ boundedness of the Riesz transform 
for $p>2$ (see \cite{CD2} and \cite{CD3}).
In Theorem~\ref{pw2} we describe a  natural generalization of the main result 
from \cite{CD2} and \cite[Theorem~5.5]{CD3}.
We obtain Theorem~\ref{pw2}  as a straightforward consequence of  
Theorem~\ref{gl} and Theorem~\ref{fsp}. 
The main idea of the proof of Theorem~\ref{gl2} comes from \cite{Si2}.
However, here we significantly simplify the proof. We also state the result
in a substantially more general and natural setting.

\section{Assumptions and notation}

Before we state our main results we have to introduce some notation 
and describe our set of assumptions.
  
\begin{assum}\label{1a}
  Let $X$ be a metric measurable space equipped with a Borel measure $\mu$ 
  and metric $\rho$. Next let $B(x,r)=\{y\in X,\, \rho(x,y)\le r\}$ be the
  open ball with center at $x$ and radius $r$.
  We suppose  throughout  that $X$ satisfies  the doubling
  property, i.e. there  exists  a constant $C$ such that
    \begin{equation}\label{doubling1.1}
    \mu(B(x, 2r)) \le C \mua 
    \end{equation}
   uniformly for all $x\in X$ and for all $r>0$.
\end{assum}
Note also that   (\ref{doubling1.1})  implies  that there   exist  positive
constants $C$ and $D$ such that
    \begin{equation}\label{d}
    \mu(B(x, \gamma r)) \le C (1+\gamma)^D \mua  \ \ \forall \gamma > 0,
    x \in X, r > 0.
    \end{equation}
In the sequel the value $D$ always refers to the constant in (\ref{d}).

 Next suppose that $TX$ is a continuous vector bundle with 
the base $X$, fibers $\C^l$
and with measurable (with respect to $x$) scalar product $(\cdo, \cdo)_x$. For 
$f(x)\in T_xX$ we put $|f(x)|_x^2=(f(x), f(x))_x$. To simplify the notation
 we will write
$(\cdo, \cdo)$   and $|\cdo|$ instead of $(\cdo, \cdo)_x$ and $|\cdo|_x$.
Now for sections $f$ and $g$ of $TX$ we put 
$$
\|f\|_{L^p(TX)}^p=\int_X|f(x)|^p\dd\mu(x) \qquad \mbox{and} \quad  \langle f,g \rangle= \int_X (f(x),g(x))\dd\mu(x).
$$
By $L^p(TX, \mu)$ we denote the Banach spaces corresponding to this norms. $L^2(TX, \mu)$ is a Hilbert space with 
the scalar product $\langle \cdo, \cdo \rangle$. 

Now suppose that $L$ is a self-adjoint positive definite operator
acting on $L^2(TX,\mu)$.
Such an operator admits a 
spectral  decomposition  $E_L(\lambda)$  and for  any  bounded  Borel
function $F\colon [0, \infty) \to \C$, we define the operator $F(L)\colon
\, L^2(TX) \to L^2(TX)$
by the formula
  \begin{equation}\label{equw}
  F(L)=\int_0^{\infty}F(\lambda) \dd E_L(\lambda).
  \end{equation}
Suppose that $S$ is a bounded operator from $L^p(TX)$ to~$L^q(TX)$.
We write $\| S \|_{L^p(TX) \to L^q(TX)}$ for the usual operator
norm of~$S$.  If $S$ is of weak type~$(1,1)$, i.e., if
$$
\mu (\{ x \in X : | Sf(x) | > \lambda \})
\leq C \, \frac{\| f \|_{L^1(TX)}}{\lambda}
\qquad\forall \lambda \in \R^+ \qquad\forall f \in L^1(TX),
$$
then we write $\| S \|_{L^1 \to L^{1,\infty}}$ for the
least possible value of $C$ in the above inequality;  this is
often called the ``operator norm'', though in fact it is not a norm.

Now let us describe the notion of \emph{integral operators}. 
For any point  $(x,y) \in X^2$ we consider the space $T_y^*\otimes T_x$.
The space $T_y^*\otimes T_x$ is canonically isomorphic to 
${\mbox {\it Hom}}\,(T_y,T_x)$, the space of all
linear homeomorphisms from $T_y$ to $T_x$. 
$T_x$ and $T_y$ are equipped with scalar product and
one can consider two natural norms on $T_y^*\otimes T_x$. These norms are: the 
Hilbert-Schmidt norm $|\cdo|_{HS}$ and the operator norm $|\cdo|$.
Let us recall that if $K(x,y)\in {\mbox {\it Hom}}\,(T_y,T_x)\cong
T_y^*\otimes T_x$, then $|K(x,y)|^2_{HS}=\mbox{\rm Tr}\,K(x,y)K(x,y)^*=
\sum_{i=1}^l\sum_{j=1}^l|(K(x,y)e_i^y,e_j^x)|^2$, where $K(x,y)^*$ is 
the adjoint of $K(x,y)$ and $e_1^z,\ldots,e_l^z$ are  arbitrary
orthonormal bases for $T_z$, $z=x$ or $z=y$.
Note that 
\begin{equation}\label{abc}
|K(x,y)|\le  |K(x,y)|_{HS}
\le l^{1/2}|K(x,y)|
\end{equation}
for all $K(x,y)\in T_y^*\otimes T_x$.
By $(T^*\otimes T) X^2$ we denote the continuous bundle with the base space 
equal to $X^2$ and with the fiber over the point $(x,y)$ equal to 
$T_y^*\otimes T_x$.
If there is a section~$K_{S}$ of  $(T^*\otimes T) X^2$
such that $|K_{S}|$ is a locally integrable function
on $(X^2,\mu\times\mu)$ and 
\begin{equation*}
  \langle Sf_1,f_2\rangle
= \int_{X} (Sf_1 ,f_2) \dd \mu
= \int_{X} (K_{S}(x,y) \, f_1(y),f_2(x)) \dd \mu(y) \dd \mu(x)
\end{equation*}
for all sections $f_1$ and $f_2$ in~$C_c(TX)$, then we say that $S$ is an
\emph{integral operator} with kernel $K_{S}$. Note that if for some $q \in 
[1,\infty)$ and for every
$x\in X$ there exists a constant $C_x$ such that 
\footnote{We assume that $C_x$ is locally integrable as a function of $x$.}
\begin{equation}\label{ke}
\sup_{f\neq 0} \frac{|Sf(x)|}{\|f\|_{L^q(TX)}}=C_x < \infty \, ,
\end{equation}
then by Riesz representation theorem (see \cite[Theorem~1, p. 286]{DuS})
 $S$ is an integral operator and 
\begin{eqnarray}\label{ke1}
l^{-1}\||K_{S}(x,\cdo)|\|_{L^{q'}(X)}\le
    l^{-1}\||K_{S}(x,\cdo)|_{HS}\|_{L^{q'}(X)} \le C_x \\ \le  
\| |K_{S}(x,\cdo)|\|_{L^{q'}(X)}
\le \||K_{S}(x,\cdo)|_{HS}\|_{L^{q'}(X)}, \nonumber
\end{eqnarray}
where $1/q+1/q'=1$ and $1\le q < \infty$. 
Thus if $S_1$ satisfies (\ref{ke}) for $q=2$ and $S_2$ is a bounded operator
on $L^2(TX,\mu)$ then $S_1S_2$ is an integral operator. Moreover,
\begin{equation}\label{ke22}
\| |K_{S_1S_2}(x,\cdo)|_{HS}|\|_{L^2(X)}\le 
\|S_2\|_{L^2(TX) \to L^2(TX)}\| |K_{S_1}(x,\cdo)|_{HS}\|_{L^2(X)}.
\end{equation}
In particular if $F_1(L)$ and $F_2(L)$  are the operators defined by 
(\ref{equw}), 
then 
\begin{equation}\label{ke2}
\||K_{F_1F_2(L)}(x,\cdo)|_{HS}\|_{L^2(X)}\le 
\|F_1\|_{L^\infty}\||K_{F_2(L)}(x,\cdo)|_{HS}\|_{L^2(X)}.
\end{equation}
Note also that 
\begin{eqnarray}\label{ke3}
\| |K_{F(L)}(x,\cdo)|\|_{L^2(X)}=\| |K_{\bar{F}(L)}(\cdo,x)|\|_{L^2(X)}, \\
\| |K_{F(L)}(x,\cdo)|_{HS}\|_{L^2(X)}=
\| |K_{\bar{F}(L)}(\cdo,x)|_{HS}\|_{L^2(X)} \quad \mbox{and} \\ 
|K_{F_1F_2(L)}(x,y)|\le \| |K_{F_1(L)}(x,\cdo)|\|_{L^2(X)}
\|| K_{F_2(L)}(\cdo,y)| \|_{L^2(X)}.\label{kke}
\end{eqnarray}
Next
\begin{equation}\label{hk1}
\mbox{\rm Tr}\,K_{|F|^2(L)}(x,x)= 
\||K_{F(L)}(x,\cdo)|_{HS}\|^2_{L^2(X)}=\||K_{\bar{F}(L)}(\cdo,x)|_{HS}\|^2_{L^2(X)}
\end{equation}
and so
\begin{equation}\label{hk2}
\mbox{\rm Tr}\,K_{\exp(-2tL)}(x,x)= 
\||K_{\exp(-tL)}(x,\cdo)|_{HS}\|^2_{L^2(X)}
=\||K_{\exp(-tL)}(\cdo,x)|_{HS}\|^2_{L^2(X)}.
\end{equation}
Finally note that by (\ref{ke1})
\begin{equation}\label{lpe1}
 l^{-1} \sup_{x\in X} \| |K_{S}(x,\cdo)|\|_{L^{q'}(X)}\le
\| S\|_{L^q(X) \to L^\infty(X)}
\le \sup_{x\in X} \| |K_{S}(x,\cdo)|\|_{L^{q'}(X)}.
\end{equation}
Hence (see \cite[Theorem~6, p. 503]{DuS})
if $S$ is bounded from $L^1(TX)$ to~$L^q(TX)$, where $q > 1$, then
$S$ is an integral operator, and
\begin{equation}\label{lpe}
l^{-1} \sup_{y\in X} \| |K_{S}(\cdo,y)|\|_{L^q(X)}\le 
\| S\|_{L^1(TX) \to L^q(TX)}
\le \sup_{y\in X} \| |K_{S}(\cdo,y)|\|_{L^q(X)};
\end{equation}
\emph{vice versa}, if $S$ is an integral operator and the right hand
side of the above inequality is finite, then $S$ is bounded from
$L^1(TX)$ to~$L^q(TX)$, even if $q = 1$.

\begin{theorem}\label{ellip} 
  Let $X$ be a measurable metric space with the doubling condition and let 
$L$ be a self-adjoint positive definite operator. 
The following conditions are equivalent:  
\begin{equation}\label{ell1} 
    \||K_{\exp(-tL)}(x,\cdo)|  \|^2_{L^2(X)} \le C\mu(B(x, t^{1/2}))^{-1} 
    \quad \forall t > 0, x \in X ;
    \end{equation}
\begin{equation}\label{ell2} 
    \| |K_{(I+tL)^{-m/4}}(x,\cdo)|  \|^2_{L^2(X)} \le C_m\mu(B(x, t^{1/2}))^{-1} 
    \quad \forall t > 0, x \in X 
    \end{equation}
for any $m>D$, where $D$ is the constant from condition~(\ref{d}) .
\end{theorem}
\begin{proof} Note that 
$$ 
{(I+(tL))^{-m/4}}
= \frac{1}{\Gamma(m/4)}\int_0^\infty e^{-s} \, s^{m/4-1}\exp(-s(tL))  
\wrt s .
$$ Hence by (\ref{d})
\begin{eqnarray*}
\| |K_{(I+tL)^{-m/4}}(x,\cdo)|  \|_{L^2(X)}
&\le& \frac{1}{\Gamma(m/4)}
\int_0^\infty e^{-s} \, s^{m/4-1}\||K_{\exp(-ts L)}(x,\cdo)|\|_{L^2(X)}
  \wrt s \\ &\le&
\frac{1}{\Gamma(m/4)}
  \int_0^\infty e^{-s} \, s^{m/4-1}\mu(B(x, (st)^{1/2}))^{-1/2} 
  \wrt s \\
&\le&  \frac{1}{\Gamma(m/4)}\mu(B(x, t^{1/2}))^{-1/2} 
\int_0^\infty e^{-s} \, s^{m/4-1}(1+1/s)^{D/4}\dd s \\&=& 
C\mu(B(x, t^{1/2}))^{-1/2}.
\end{eqnarray*}
To prove that (\ref{ell2}) implies (\ref{ell1}) we note that by (\ref{ke2})
and (\ref{abc})
\begin{eqnarray}\label{ondi}
\||K_{\exp(-t L)}(x,\cdo)|\|_{L^2(X)}
\le l^{1/2}\|\exp(-tL)(1+tL)^{m}\|_{L^2 \to L^2}
\| |K_{(I+tL)^{-m}}(x,\cdo)|  \|_{L^2(X)}\nonumber \\ \hspace{1.8cm}
\le l^{1/2} \sup_{\lambda \in \R^+} e^{-t\lambda}(1+t\lambda)^m 
\| |K_{(I+tL)^{-m}}(x,\cdo)|  \|_{L^2(X)}\le C 
\| |K_{(I+tL)^{-m}}(x,\cdo) | \|_{L^2(X)}.
\end{eqnarray}\end{proof}

{\em Remarks.} 1. Note that Theorem~\ref{ellip} remains valid if 
we replace $\mu(B(x, t^{1/2}))^{-1/2}$ by $v_x(t)$ for any decreasing function $v_x$. 
(\ref{ell2}) implies (\ref{ell1}) without any assumptions on 
$v_x$ or $\mu$. To show the inverse implication
one has to assume that $v_x(ts)\le C v_x(t)(1+1/s)^{D/4}$. For 
$v_x(t)=  \mu(B(x, t^{1/2}))^{-1/2} $ this means that $\mu$ satisfies
condition (\ref{d}).

2. Note that in virtue of (\ref{hk2}) and (\ref{abc})
it is enough to know  the value of $\mbox{\rm Tr}\,K_{\exp(-2tL)}(x,x)$ to 
verify condition~(\ref{ell1}). Therefore
condition~(\ref{ell1}) is often called on-diagonal bounds of a heat kernel.
Condition~(\ref{ell1}) is well understood.
On-diagonal bounds are very often used as
a basic assumption in theorems concerning the heat kernels and
boundedness of the Riesz transforms (see for example \cite{CD1, DP, Gr1}).
There are many examples of operators
satisfying  condition~(\ref{ell1})  and there 
are efficient techniques to obtain 
condition~(\ref{ell1})
for some particular class of operators in the scalar case. For
 example it is known that
for the Laplace-Beltrami 
operators condition~(\ref{ell1}) is equivalent to a relative Faber-Krahn
inequality (see \cite{Gr}).
The related literature is too large to be listed here, so we refer reader to
\cite{Co2, Da, Gr, Ro, VSC} for the related theory, examples of 
operators satisfying (\ref{ell1}) and for further references.

\medskip
Now we set 
\begin{equation*}
\D_r=\{ (x,y)\in X\times X: \rho(x,y) \le r \}.
\end{equation*}
Given an operator $S$ from $L^p(TX)$ to $L^q(TX)$, we write
\begin{equation}\label{kernel}
\supp K_{S} \subseteq \D_r
\end{equation}
if $\langle S f_1, f_2 \rangle = 0$ whenever $f_n$ is in~$C(TX)$
and $\supp f_n \subseteq B(x_n,r_n)$ when $n = 1,2$, and
$r_1+r_2+r
< \rho(x_1, x_2)$.  This definition makes sense even if $S$ is
not an integral operator, in the sense of the previous definition.
If S is an integral operator with the kernel $K_S$, 
then (\ref{kernel}) is equivalent to the standard meaning of 
$\supp K_{S}  \subseteq \D_r$,
 that is $K_S(x,y)=0$ for all $(x,y) \notin \D_r$.

\begin{theorem}\label{fspro}
 Let $L$ be a self-adjoint positive definite operator acting on $L^2(X)$. 
The following conditions are equivalent: 
\begin{equation}\label{fs1}
\supp K_{\cos(t\sqrt L)} \subseteq \D_t \quad \forall t\ge 0\,;
\end{equation}
\begin{equation} \label{fs2}
   | \langle \exp(-tL)f_1,f_2\rangle | \le Ce^{-\frac{r^2}{4t}}
\|f_1\|_{L^2(TX)} \|f_2\|_{L^2(TX)}
\quad \forall
t >0, 
    \end{equation}
whenever $f_n$ is in~$C(TX)$
and $\supp f_n \subseteq B(x_n,r_n)$ when $n = 1,2$, and
$0\le r  
< \rho(x_1, x_2)- (r_1+r_2)$.
\end{theorem}
{\em Remark.} The connection of the heat and the wave equation has a long
history (see \cite{CGT, Me}, see also \cite{Si2}
 and the third proof of \cite[Theorem~3.2, p. 157]{Gr}). For the origin of
the $L^2$ Gaussian estimates (\ref{fs2}) so-called the Davies 
or the Davies-Gaffney
estimates see \cite{Da11}.

\begin{proof}
Suppose that $\supp f_n \subseteq B(x_n,r_n)$ for $n = 1,2$, and that
$0\le r < \rho(x_1, x_2)- (r_1+r_2)$. 
Put 
$$
u(z)=<\exp(-{L}/{(4z)})f_1,f_2>.
$$ 
$L$ is a self-adjoint positive definite operator so
$u$ is an analytic function on the complex 
 half-plane $\mbox{Re}\,z > 0$, continuous and bounded
on the set  $\{z\in \C \colon \mbox{Re}\,z \ge 0, \, z\neq 0\}$, and
\begin{equation*}
\sup {}_{\re z=0}\,|e^{r^2z}u(z)|\le 
\|f_1\|_{L^2(TX)} \|f_2\|_{L^2(TX)}.
\end{equation*}
By (\ref{fs2}) 
$$
\sup{}_{z \in \R_+}\,
|e^{r^2z}u(z)| \le C\|f_1\|_{L^2(TX)} \|f_2\|_{L^2(TX)}.$$
 Hence, by Phragm\'en-Lindel\"of theorem for an angle
(see \cite[Theorem~7.5, p.214, vol. II]{Ma} or \cite[Lemma~4.2, p.107]{StW}) 
$$|e^{r^2z}u(z)|  \le  
\|f_1\|_{L^2(TX)} \|f_2\|_{L^2(TX)} $$
and
\begin{equation}\label{anly}
|u(z)|\le \exp(-r^2\mbox{Re}\,z) \|f_1\|_{L^2(TX)} \|f_2\|_{L^2(TX)}
\end{equation}
for all $z$ such that $\mbox{Re}\,z > 0$.
Next we note that
\begin{equation}\label{fc}
<\exp(-s L)f_1,f_2>= \int_0^{\infty} < \cos(t\sqrt L)f_1,f_2>
\frac{2}{\sqrt {\pi s}} e^{-\frac{t^2}{4s}} \dd t.
\end{equation} 
By change of  variable $t:=\sqrt t$ in integral (\ref{fc}) and putting 
$s:=1/(4s)$ we get 
\begin{equation}\label{dok}
{s}^{-1/2}<\exp{\Big(-\frac{L}{4s}\Big)}f_1,f_2>= 2\int_0^{\infty} 
({\pi t})^{-1/2} <\cos(\sqrt t \sqrt L)f_1, f_2>
 e^{-st} \dd t,
\end{equation}
so the function  $v(z)={z}^{-1/2}u(z) $ is a 
Fourier-Laplace transform of the function
$w(t)=(\sqrt{\pi t})^{-1} < \cos(\sqrt t \sqrt L)f_1, f_2>$.
Now by (\ref{anly}) and  the Paley-Wiener Theorem 
(Theorem 7.4.3 \cite{Ho})
\begin{equation}\label{koniec}
{\rm supp}\; w \subseteq [r^2, \infty).
\end{equation}
This proves that (\ref{fs2}) implies (\ref{fs1}).
Now if (\ref{fs1}) holds, then by (\ref{fc})
\begin{eqnarray*}
|< \exp(-s L)f_1,f_2>| \le 
\int_0^{\infty} |< \cos(t\sqrt L)f_1,f_2>|
\frac{2}{\sqrt {\pi s}} e^{-\frac{t^2}{4s}} \dd t  \hspace{3cm}{}\\=
\int_r^{\infty} |< \cos(t\sqrt L)f_1,f_2>|
\frac{2}{\sqrt {\pi s}} e^{-\frac{t^2}{4s}} \dd t\le
 \|f_1\|_{\ld}\|f_2\|_{\ld}\int_r^{\infty}
\frac{2}{\sqrt {\pi s}} e^{-\frac{t^2}{4s}} \dd t \\ \le
e^{-\frac{r^2}{4s}}\|f_1\|_{\ld}\|f_2\|_{\ld}.
\end{eqnarray*}
\end{proof} 

The following lemma is a very simple
but useful consequence of (\ref{fs1}).
\begin{lemma}\label{step}
Assume that $L$ satisfies {\rm (\ref{fs1})} and that
 $\widehat{F}$ is the Fourier transform
of an  even bounded Borel function $F$ with 
$\mbox{\rm supp}\; \widehat{F} \subseteq [-r,r]$.
Then 
$\mbox{\rm supp}\; K_{F(\sqrt L)} \subseteq \D_r$.
\end{lemma}
\begin{proof}
Since $F$ is even, by the Fourier inversion formula,
$$
F(\sqrt L) =\frac{1}{2\pi}\int_{-\infty}^{+\infty}%
  \widehat{F}(t) \cos(t\sqrt L) \;dt.
$$
But  supp $\widehat{F} \subseteq [-r,r]$ 
and Lemma \ref{step} follows from (\ref{fs1}).
\end{proof}

\section{Main results}
We are now in a position to state our two main results.

\begin{theorem}\label{gl2}
Suppose that for some number $N\in \N$ and points $x,y \in X$
there exist functions $V_x,V_y \colon \R^+ \mapsto \R$ such that
\begin{equation}\label{aszw1}
    \| |K_{(I+t^2L)^{-N/4}}(z,\cdo)|  \|_{L^2(X)} \le  V_{z}(t)
    \quad \forall t > 0, \,\, z=x,y.
    \end{equation}
Then, there exists a constant $C_N$ such that for all $t<\rho(x,y)^2$
\begin{eqnarray}\label{zw1}
|K_{\exp(-tL)}(x,y)| \le C_N  
V_x\bigg(\frac{t}{\d(x,y)}\bigg)V_y\bigg(\frac{t}{\d(x,y)}\bigg)
 \frac{{\exp}\left(\frac{-\d(x,y)^2}{4t}\right)}{\d(x,y)t^{-1/2}}.
\end{eqnarray}
Thus if $L$ satisfies (\ref{ell1}) or (\ref{ell2}), then
\begin{eqnarray}\label{zw2}
\qquad |K_{\exp(-tL)}(x,y)| \le C  
\mu \Big(B \Big(x,\frac{t}{\d(x,y)}\Big)\Big)^{-\frac{1}{2}}
\mu \Big(B \Big(y,\frac{t}{\d(x,y)}\Big)\Big)^{-\frac{1}{2}}
 \frac{\exp{\Big(\frac{-\d(x,y)^2}{4t}\Big)}}{\d(x,y)t^{-1/2}}
\end{eqnarray}
for all $t<\rho(x,y)^2$.
\end{theorem}

\begin{theorem}\label{gl}
Suppose that $X$ is a measurable metric space satisfying
Assumption~\ref{1a} and that $TX$ and $T'X$ are vector bundles with 
measurable scalar products. Next assume that an  operator  $L$ acting
on $L^2(TX)$ satisfies {\rm(\ref{ell2})}
 and {\rm
(\ref{fs1})}. 
Assume also that $A \colon \, D(A) \to L^2(T'X)$ is a local operator, which
means that for any $f\in D(A)\subset L^2(TX)$ 
\begin{equation}\label{lokal}
\supp Af \subseteq \supp f.
\end{equation}
Finally assume that $D(L^{\alpha}) \subset D(A)$ and 
$AL^{-\alpha}\colon \, L^2(TX) \to L^2(T'X) $ is
bounded  for some $\alpha >0$.
Then the operator $AL^{-\alpha}$ is of weak type $(1,1)$ and bounded 
as an operator from  $L^p(TX,\mu)$ to $L^p(T'X,\mu)$
for all $p\in (1,2]$.
\end{theorem}

{\em Remarks.} 1. Note that by (\ref{kke}) and (\ref{ondi})
$
|K_{\exp(-2tL)}(x,y)|\le V_x(\sqrt t)V_y(\sqrt t).
$
For   $t \ge \rho(x,y)^2$ this obvious 
estimate  is sharp even for the standard Laplace operator. 
Therefore the discussion of Gaussian bounds
for $t \ge \rho(x,y)^2$ is straightforward so we do not include description
of the heat kernel bounds  for this case 
in  the statement of Theorem~\ref{gl2}.

2. Theorem~\ref{gl2} holds without the doubling volume property. One needs 
 Assumption~\ref{1a} only to prove that (\ref{ell1}) and (\ref{ell2})
are equivalent. In our proof that (\ref{aszw1}) implies (\ref{zw1}) 
and that (\ref{zw2}) follows from (\ref{ell1}) we do not require any
condition on $V_z$ or $\mu(B(x,t))$. Note also that if
$V_z$ satisfies the doubling condition, then estimates (\ref{aszw1})
are equivalent to the on-diagonal estimates for the heat kernel 
$ \| |K_{\exp{(t^2L)}}(z,\cdo)|  \|_{L^2(X)} \le C V_{z}(t)$ 
(see {\em Remark}~1. Theorem~\ref{ellip}).

3. Note that if Assumption~\ref{1a} holds
then 
$$\mu(B(x,{t}/{\d(x,y)}))^{-\frac{1}{2}}\le 
\mu(B(x,{\sqrt t}))^{-\frac{1}{2}}(\d(x,y)/\sqrt t)^{D/2}$$ and by (\ref{zw2})
\begin{eqnarray}\label{mol1}
|K_{\exp(-tL)}(x,y)| \le C  
\mu(B(x,{\sqrt t})^{-\frac{1}{2}}
\mu(B(y,{\sqrt t}))^{-\frac{1}{2}}{(\d(x,y)/\sqrt{t})^{(D-1)/2}}
{\exp{\Big(\frac{-\d(x,y)^2}{4t}\Big)}}.
\nonumber
\end{eqnarray}
In \cite{Mo}  Molchanov  proved 
that if N is the north pole and S is the south  pole of the $D$-dimensional
 unit sphere and $L$ is the Laplace-Beltrami operator on the sphere, then
\begin{equation}\label{mol2}
K_{exp(-tL)}(N,S) \sim  t^{-D/2}(1+ \rho(S,N)/\sqrt t)^{D-1}
 {\exp}\left(-\frac{\rho(N,S)^2}{4t}\right) \;\;\;\; {\rm as} \;\; t \downarrow 0.
\end{equation}
This shows that  estimates  (\ref{zw1}) and (\ref{zw2}) are sharp (see also
\cite[Theorem 5.9]{Gr} and \cite{Si2}).

\section{Finite speed propagation}

Finite speed propagation property for the solution of the wave equation is
one of our main assumptions. Hence for the sake of completeness we describe the
proof of finite speed propagation property for a large class of operators.

Finite speed propagation property for the wave propagator is well known
(see for example \cite[Theorem~(5.3)]{Fo}, \cite[Theorem~6.1]{Ta}).  
However, the statement of Theorems~\ref{hod1} below is more general
than what is usually found in the literature. Moreover, the proof given
 here is simpler than other proofs known to the 
author.

Suppose that $M$ is a complete Riemannian manifold and $\mu$ is an absolutely
continuous measure with a smooth density not equal to zero at any point of $M$.
By $\Lambda^kT^*M$ we denote the bundle of $k$-forms on $M$.
For fixed  $\beta, \beta_* \in L^2(\Lambda^1T^*M)$ 
and $\gamma \in L^2(\Lambda^kT^*M)$
We define the operator $L$ ($L=L_{\beta,\beta_*,\gamma}$) acting on 
$L^2(\Lambda^kT^*M$) by the formula
\begin{equation}\label{hod1}
\langle L \omega, \omega\rangle =
\int_M |d_k\omega +\omega \wedge \beta|^2+
|d_{n-k}*\omega +*\omega \wedge \beta_*|^2+
|*\omega \wedge \gamma|^2
\dd\mu(x),
\end{equation}
where $\omega$ is a smooth compactly supported $k$-form
and $*$ is the Hodge star operator. With some abuse of notation we also 
denote by $L$ its Friedrichs extension.
Note that for example the Hodge Laplacian (see \S~\ref{rhl}) and 
Schr\"odinger operators with  electromagnetic fields (see \S~\ref{sh})
can be defined by (\ref{hod1})).

\begin{theorem}\label{fsp}
The operator $L$ defined by (\ref{hod1}) satisfies (\ref{fs1})
and (\ref{fs2}).
\end{theorem}
\begin{proof}
We put $\omega_t(x)=\omega(t,x)=\exp(-tL)\omega$. Then we fix some function 
$\xi\in C^\infty(M)$ such that $|d \xi|\le \kappa$ 
and we consider the integral
$$
E(t)=\int_M(\omega(t,x),\omega(t,x))e^{\xi(x)}\dd\mu(x).
$$
Next we note that for every $k$-form $\eta$ and $1$-form $\zeta$ we have
$|\zeta \wedge \eta|^2+|\zeta \wedge *\eta|^2=|\eta|^2|\zeta|^2$ and
\begin{eqnarray*}
\frac{E'(t)}{2}=\re \int_M(\partial_t\omega(t,x),\omega(t,x))e^\xi\dd\mu(x)
=-\re \int_M(L\omega_t,\omega_t e^\xi)\dd\mu\\
=-\re \int_M \big[(d_k\omega_t+\omega_t\wedge \beta,%
 d_k (\omega_t e^\xi) +(\omega_t e^\xi)  \wedge \beta)
+(*\omega_t\wedge \gamma, *\omega_t\wedge \gamma)e^\xi\big]
\dd\mu
\\-\re \int_M
(d_{n-k}*\omega_t +*\omega_t\wedge \beta_*, d_{n-k} (*\omega_t e^\xi)+ 
(*\omega_t e^\xi)\wedge \beta_*)\dd\mu\\
=-\int_M \big[|d_k\omega_t +\omega_t\wedge \beta|^2 e^\xi+
 |d_{n-k}*\omega_t+*\omega_t\wedge \beta|^2 e^\xi+
|*\omega_t\wedge \gamma|^2e^\xi\big]
\dd\mu\\
-\re \int_M (d_k\omega_t+\omega_t\wedge \beta, d_0\xi \wedge \omega_t) e^\xi\dd\mu-
\int_M (d_{n-k}*\omega_t+*\omega_t\wedge \beta_*,d_0\xi \wedge  *\omega_t) e^\xi\dd\mu
\\ \le
-\int_M \big[|d_k\omega_t +\omega_t\wedge \beta|^2 e^\xi+
 |d_{n-k}*\omega_t+*\omega_t\wedge \beta_*|^2 e^\xi+
|*\omega_t\wedge \gamma|^2e^\xi\big]
\dd\mu
\\
+\int_M \big[|d_k\omega_t +\omega_t\wedge \beta|^2 e^\xi+
 |d_{n-k}*\omega_t+*\omega_t\wedge \beta_*|^2 e^\xi+
|*\omega_t\wedge \gamma|^2e^\xi\big]
\dd\mu
\\  
+ \frac{1}{4}\int_M \big[|d_0\xi \wedge \omega_t |^2 e^\xi 
+ |d_0\xi \wedge *\omega_t |^2 e^\xi\big]\dd\mu
= \frac{1}{4}\int_M |\omega_t|^2|d_0\xi|^2 e^\xi\dd\mu\le 
\frac{\kappa^2 E(t)}{4}\,.
\end{eqnarray*}
Hence $E(t)\le \exp(\kappa^2t/2)E(0)$.
Now we say that $\xi\in \Theta_\kappa \subseteq C^{\infty}(M)$ 
if $\xi(x)=0$ for $x\in B(x_1,r_1)$ and $|d\xi| \le \kappa$.
Next assume  that
$0\le r < \rho(x_1, x_2)- (r_1+r_2)$. Then 
$$
\sup_{\xi\in \Theta_\kappa}
\int_{B(x_2,r_2)}|\omega|^2 e^{\xi}\dd\mu\ge e^{r\kappa}
\int_{B(x_2,r_2)}|\omega|^2\dd\mu.
$$
Hence if supp~$\omega_0 \subseteq B(x_1,r_1)$
then
$$\int_{B(x_2,r_2)}|\omega_t|^2\dd\mu
\le \exp\Big(\frac{\kappa^2t}{2}-r\kappa\Big)
\int_{M}|\omega_0|^2\dd\mu 
$$ 
Putting $\kappa=r/t$ we get
\begin{equation}\label{kar}
\int_{B(x_2,r_2)}|\omega_t|^2\dd\mu\le \exp\Big(\frac{-r^2}{2t}\Big)
  \int_{M}|\omega_0|^2\dd\mu  = \exp\Big(\frac{-r^2}{2t}\Big)
\int_{B(x_1,r_1)}|\omega_0|^2\dd\mu
\end{equation}
Now (\ref{fs2}) is a straightforward consequence of (\ref{kar}).
\end{proof}

\section{Off-diagonal Gaussian bounds, proof of Theorem~\ref{gl2}}

\begin{proof}
For $s>1 $, we define the family of functions $\phi_s$ by the formula
$$
\phi_s(x)= \psi(s(|x|-s)),
$$
where $\psi \in C^{\infty}(\R)$ and 
$$
\psi(x) = \left\{ \begin{array}{ll}
    0 & \mbox{ if $x \le -1$}\\
    1  & \mbox{ if $x \ge -1/2$}\; .
           \end{array}
    \right.
$$
Finally we define functions $F_s$ and $R_s$ by the following formula
\begin{equation*}
F_s(x) 
= \frac{1}{\sqrt{4\pi}}\exp {(\frac{-x^2}{4})}-R_s(x)
= \phi_s(x) \frac{1}{\sqrt{4\pi}}\exp {(\frac{-x^2}{4})}
\end{equation*}
so that $\widehat{F_s}(\lambda)+\widehat{R_s}(\lambda)=\exp(-\lambda^2)$
and
\begin{equation}\label{ggg}
\widehat{F_s}(\sqrt{tL})+\widehat{R_s}(\sqrt{tL})=\exp(-tL).
\end{equation}
Integration by parts $N$ times yields
\begin{eqnarray*}
\int \phi_s(x)e^{-\frac{x^2}{4}}e^{-ix\lambda}=
\int\underbrace{\Big(\frac{1}{x/2+i\lambda}%
\Big(\ldots\Big(\frac{1}{x/2+i\lambda}\phi_s(x)\Big)^{\prime}\ldots\Big)'\Big)%
^{\prime}}_Ne^{-\frac{x^2}{4}-i\lambda x}dx. 
\end{eqnarray*}
Hence for any natural number $N$ and $s>1$
\begin{equation}\label{osz}
|\widehat{F_s}(\lambda)| \le 
C'_N\frac{1}{s(1+\lambda^2/s^2)^{N/2}}e^{-\frac{s^2}{4}} ,
\end{equation}
where $C'_N$ is a constant depending only on $N$.
Next we note that 
supp~$R_s \subseteq [-s+\frac{1}{2s},s-\frac{1}{2s}]$,
so if we put $s_{xy}= \rho(x,y)t^{-1/2} $, then 
Lemma~\ref{step}
$K_{\widehat{R_{s_{xy}}}(\sqrt {tL})}(x,y)=0$.  Hence by (\ref{ggg})
\footnote{(\ref{incl}) shows that the remainder 
${\widehat{R_{s_{xy}}}(\sqrt {tL})}$ does not contribute to the value of
the heat kernel $K_{\exp(-tL)}(x,y)$. Subtracting the remainder from the 
heat propagator is the main idea of the proof.}
\begin{equation}\label{incl}
K_{\exp(-tL)}(x,y)=K_{\widehat{F_{s_{xy}}}(\sqrt {tL})}(x,y).\;\; 
\end{equation}
Now let $J_{s_{xy}}$ be a function such that  
$J_{s_{xy}}(\lambda)^2=\widehat{F_{s_{xy}}}(t^{1/2}\lambda)$.
By (\ref{osz})
$$
\sup_{\lambda \ge 0}\bigg|J_{s_{xy}}(\lambda)
\bigg(1+\frac{\lambda^2t^2}{\d(x,y)^2}\bigg)^{N/4}\bigg|\le 
C
\frac{\exp\left(-\frac{\d(x,y)^2}{8t}\right)}{\sqrt{\d(x,y)t^{-1/2}}}\,.
$$
Hence by (\ref{ke2}) and (\ref{abc})
\vspace{-0.1cm}
\begin{eqnarray}
 \| |K_{J_{s_{xy}}(\sqrt L)}(x,\cdo)|  \|_{L^2(X)} &\le& C 
\frac{\exp\left(-\frac{\d(x,y)^2}{8t}\right) 
}{\sqrt{\d(x,y)t^{-1/2}}}
\Big\|\Big|K_{\left(I+\frac{t^2L}{\d(x,y)^2}\right)^{-N/4}}(x,\cdo)
\Big|\Big\|_{L^2(X)} \nonumber
\\ &\le& C 
V_x\bigg(\frac{t}{\d(x,y)}\bigg)
\frac{\exp\left(-\frac{\d(x,y)^2}{8t}\right)}{\sqrt{\d(x,y)t^{-1/2}}}\,.
\label{hcaa}
\end{eqnarray}
Finally by (\ref{kke})
\begin{equation}\label{hcab}
|K_{\exp(-tL)}(x,y)|=|K_{\widehat{F_{s_{xy}}}(\sqrt {tL})}(x,y)|
\le \| K_{J_{s_{xy}}(\sqrt L)}(x,\cdo)  \|_{L^2(X)}
\| K_{J_{s_{xy}}(\sqrt L)}(y,\cdo)  \|_{L^2(X)}
\end{equation}
and (\ref{zw1}) follows from (\ref{hcaa}) and (\ref{hcab}).
\end{proof}

\section{Riesz Transform, proof of Theorem~\ref{gl}}

We fix an even function~$\Phi$ in the Schwartz
space~$\mcS(\R)$ such that $\Phi(0)
= 1$, whose Fourier transform $\widehat\Phi$ is supported in~$[-1,1]$;
we let $\Phi_{r}$ denote the dilated function~$\Phi(r\cdot\,)$
and $\Phi^{(l)}$ denote the $l^{\,\mathrm{th}}$ derivative of $\Phi$.
For later purposes, note that for any fixed positive integer
$K$, one may assume that $\Phi^{(l)} (0) = 0$ when $1 \leq l \leq
K$.
\begin{lemma} \label{lns805.4.1}
Let $\Phi$ in~$\mcS(\R)$ be chosen as above.
If {\rm (\ref{ell2})} and {\rm (\ref{fs1})} hold, then the kernel
$K_{\Phi_{{r}} (\sqrt L)}$ of the self-adjoint 
operator $\Phi_{{r}} (\sqrt L)$ satisfies 
\begin{equation}
{\rm supp}\; K_{\Phi_{{r}} (\sqrt L)}
 \subseteq   \D_r   
\label{elp}
\end{equation}
and
\begin{equation}\label{98}
\int|K_{\Phi_{{r}}(\sqrt L)}(x,y)|^2\dd\mu(x) =
\int|K_{\Phi_{{r}}(\sqrt L)}(y,x)|^2\dd\mu(x)
\leq C\, \mu(B(y,r))^{-1}
\end{equation}
for all $r>0$ and $y \in X$.  
\end{lemma}
\begin{proof} (\ref{elp}) follows from Lemma~\ref{step}. 
Next by (\ref{ke22})
$$
\||K_{\Phi_{{r}} (\sqrt L)}(\cdo,y)|\|_{L^2(X)} \ \le l^{1/2}
\|(I+r^2L)^{m} \Phi_{{r}} (\sqrt L)\|_{\ld \to \ld} \, 
\|K_{(I+r^2L)^{-m}}(\cdo,y)\|_{L^2(X)}.
$$
The $L^2$ operator norm of the first term is equal to the $L^{\infty}$
norm of the function \linebreak 
$ (1 + r^2 \lambda )^{2m} \Phi (r {\sqrt{\lambda}})$
which is uniformly bounded in $r>0$ for any fixed  $m\in \N$ and so 
(\ref{98}) follows from (\ref{ell2}).
\end{proof}
We now recall  
the  Calder\'on-Zygmund decomposition.
(see e.g. \cite{mchr2, CW, St}).
\begin{theorem}\label{calzyg}
There exists $C$ such that, given $ f \in L^1(TX,\mu)$ and $\lambda >0$,
one can decompose $f$ as
$$
f=g+b=g+\sum b_i
$$
so that
\begin{enumerate}

\item $|g(x)| \le C\lambda$, a.e. $x$ and $\|g\|_{\lj} \le C \|f\|_{\lj}$.

\item There exists a sequence of balls $B_i=B(x_i,r_i)$ such that
the support of each $b_i$ is contained in $B_i$ and
$$
\int|b_i(x)|\dd\mu(x) \ \le \ C \lambda\mu(B_i). 
$$

\item $\sum \mu(B_i) \le C\, {\frac{1}{\lambda}} \int|f(x)|\,\dd\mu(x)$.

\item  There exists $\sigma \in \N$ such that each point of $X$ is contained in
at most $\sigma$ of the balls $B(x_i, 2 r_i)$.
\end{enumerate}
\end{theorem}
The proof of Theorem~\ref{calzyg} is a variant of standard arguments, for which see, e.g. 
\cite[p. 66]{mchr2}, \cite{CW} or \cite[p. 8]{St}.\footnote{Note that we do not have to
assume that $\int b_i=0$ which could be difficult to achieve for
the vector bundle version which we consider here (see \cite{CD3}).}

\vskip 7pt
We are now in a position to prove Theorem \ref{gl}.

\begin{proof}[Proof of Theorem~\ref{gl}]
We start by decomposing $f$ into $g + \sum b_i$ at level
$\lambda$ according to Theorem \ref{calzyg}. We will follow the
idea of using more information of the $L^2$ 
operator norm (in our case, $\| AL^{-\alpha} \|_{\ld \to L^2(T'X) } <
\infty $)
by smoothing out the bad functions $b_i$ at a scale of
the size of their support and considering this part of the good function
where $L^2$ estimates can be used (see \cite{CD1, DM, F, He} for other variants
of this).

In our case  let         
$G = g + \sum \Phi_{r_i}(\sqrt{L}) b_i$ be the modified 
good function.\footnote{$\Phi_{r_i}$ is the function from 
Lemma~\ref{lns805.4.1}}
Hence $f = G + B$, where $B = \sum (I - \Phi_{r_i}(\sqrt{L})) b_i$ and
we write
\begin{eqnarray}\nonumber
\lambda \, \mu(\{ |AL^{-\alpha} f(x) | \ge \lambda \}) \ \le &
\hspace{-0.5cm}\lambda \,
\mu(\{| AL^{-\alpha} G(x) | \ge \lambda / 2 \}) \\ &+ \ \lambda \,
\mu(\{ |AL^{-\alpha}B(x) | \ge \lambda / 2\}).\label{gb}
\end{eqnarray}
The first term is less than $ 4/{\lambda} \, \|AL^{-\alpha} G \|^2_{
L^2(T'X)} \, \le \,  C/{\lambda} \, \|G\|^2_{\ld}$. However,
according to Lemma~\ref{lns805.4.1},
$\supp \Phi_{r_i}(\sqrt{L})b_i \subseteq B(x_i,2r_i)$ and by (\ref{lpe})
\begin{eqnarray*}
& \hspace{-3cm} \|\Phi_{r_i}(\sqrt{L})b_i\|_{\ld}^2 \ \le C 
\sup_{y\in B(x_i,2r_i)}\int | K_{\Phi_{r_i}(\sqrt{L})}(x,y)|^2
\dd\mu(x) 
\|b_i\|_{\lj}^2
\\  & \hspace{1cm}
\le C \mu(B(x_i,r_i))^{-1} \|b_i\|^2_{\lj} 
\ \le \
C\, \lambda \|b_i\|_{\lj}.
\end{eqnarray*}
Hence by Theorem~\ref{calzyg} again
\begin{eqnarray*}
\|G \|_{\ld}^2
\le C \Big( \|g\|_{\ld}^2 + \sigma \lambda \sum_i \|b_i\|_{\lj}\Big)
 \le C\lambda \|f\|_{\lj}.
\end{eqnarray*}
and so the first term in (\ref{gb})
is bounded by $C \, \|f\|_{\lj}$.
\vskip 5pt
Since $\mu( \cup B(x_i , 2 r_i ) ) \le C 
\sum \mu ( B_i ) \le C \|f\|_{\lj}/\lambda $, then to bound the second term in (\ref{gb}), it
suffices to show
\begin{equation}\label{fin1}
\int\limits_{x\notin \cup B^{*}_{i}} | AL^{-\alpha} B(x) |\, 
\dd\mu(x) \ \le \ C \| f \|_{\lj},
\end{equation}
where $B^{*}_{i} = B(x_i , 2 r_i)$.
The left side of (\ref{fin1}) 
is less than
\begin{eqnarray*}
& \hspace{-2cm}\sum\limits_{i} \int\limits_{x\notin \cup_{j} B^{*}_{j}} \Bigl| 
\int K_{AL^{-\alpha}(1-\Phi_{r_i})(\sqrt{L})}(x,y) b_i (y) \,
\dd\mu(y)\Bigr|\,
\dd\mu(x)\\
& \le \sum\limits_{i} \int |b_i (y) | \int\limits_{x\notin B^{*}_{i}}
| K_{AL^{-\alpha}(1-\Phi_{r_i})(\sqrt{L})} (x,y) | \, \dd\mu(x) \,
\dd\mu(y) .
\end{eqnarray*}
By (\ref{lpe})
(\ref{fin1}) follows
from Theorem \ref{calzyg} once we establish
\begin{equation}\label{fin2}
\sup\limits_{y, r} \int\limits_{\rho(x,y) \ge \,  r}
| K_{AL^{-\alpha}(1-\Phi_{r} )(\sqrt{L})} (x,y) | \, \dd\mu(x) \ \le \
C \ .
\end{equation}
By the Cauchy-Schwartz inequality, 
\begin{eqnarray}
\int\limits_{\d(x,y)\ge r}|K_{AL^{-\alpha}(1-\Phi_{r})(\sqrt L))}(x,y)| \, 
\dd\mu(x)
= 
\sum _{j\ge 1}\int\limits_{2^{j}r \ge \rho(x,y)\ge 2^{j-1} r} 
|K_{AL^{-\alpha}(1-\Phi_{r})(\sqrt L)}(x,y)|\dd\mu(x)
\nonumber\\
\le  
\sum_{j\ge 1} 
\mu(B(y,2^{j}r))^{1/2} \Big(\int\limits_{\rho(x,y)\ge 2^{j-1} r} 
|K_{AL^{-\alpha}(1-\Phi_{r})(\sqrt L)}(x,y)|^2 \dd\mu(x)\Big)^{1/2}\ .
\label{summ}
\end{eqnarray}
Fix a nonnegative even $\varphi \in C^{\infty}_c (\R )$ such that
$\varphi = 1$ on  $[-1/4,1/4]$ and $\varphi = 0$ on
$\R\setminus[-1/2,1/2]$. Set $\varphi_s(\lambda)=\varphi(s\lambda)$
and let denote by $\check{\varphi}_s$ the inverse Fourier transform 
of $\varphi_s$. We put  $H^\alpha(\lambda)=\lambda^{-2\alpha}$.
Note that $H^{\alpha} (1 - \Phi_{r})(\lambda)= 
r^{2\alpha}H^{\alpha} (1 - \Phi_{1})(r\lambda)$.
We define functions $F_j^\alpha$  and $R_j^\alpha$ by the formula
\footnote{We estimate the last integral in (\ref{summ}) in a similar way to
the proof of theorem~\ref{gl2}.}
$$
r^{2\alpha}R_j^{\alpha}(r\lambda)=[ H^{\alpha} (1 - \Phi_{r})] * 
{\check{\varphi}_{2^{-j}/r}}(\lambda)= [ H^{\alpha} (1 - \Phi_{r})]-
  r^{2\alpha}F_j^{\alpha}(r\lambda).
$$
Then 
$$ 
\supp \widehat{R_j^{\alpha}}\subset [-2^{j-1},2^{j-1}].
$$
Hence by (\ref{lokal}) and Lemma \ref{step} the kernels of 
$AH^{\alpha} (1 - \Phi_{r})(\sqrt{L})$ and 
$r^{2\alpha}AF_j^{\alpha}(r\sqrt{L})$ 
coincide on the set 
$X^2 \setminus \D_{r2^{j-1}}$ and 
\begin{eqnarray*}
\int\limits_{r2^{j-1} < \rho(x,y)} |K_{AH^{\alpha}
(1 - \Phi_{r})(\sqrt{L})} (x,y) |^2 \, \dd\mu(x) \ = \
\int\limits_{r2^{j-1} < \rho(x,y)} 
|K_{A r^{2\alpha} F_j^{\alpha}(r\sqrt{L})} (x,y) |^2 \, \dd\mu(x) \ \\ 
 \le \int_X 
|K_{A r^{2\alpha} F_j^{\alpha}(r\sqrt{L})} (x,y) |^2 \, \dd\mu(x) \  
=
\int_X |K_{AL^{-\alpha}(r^2L)^{\alpha}F_j^{\alpha}(r\sqrt{L})} (x,y) |^2 
\, \dd\mu(x).  
\end{eqnarray*} 
However, we assume that $L$ satisfies (\ref{ell2}) and by (\ref{ke22}) and (\ref{ke2})
\begin{eqnarray*}
\,\int_X \!|K_{AL^{-\alpha}(r^2L)^{\alpha}F_j^{\alpha}(r\sqrt{L})}(x,y)|^2\!  \dd\mu(x) 
\le \|AL^{-\alpha}\|^2_{L^2(TX) \to L^2(T'X)}
\| |K_{(r^2L)^{\alpha}F_j^{\alpha}(r\sqrt{L})}(\cdo,y)|\|^2_{L^2(X)}\\
\le
C \|J^{\alpha}_j(r\sqrt L) \|^2_{\ld \to \ld}
\| |K_{(I + 2^{2j}r^2L)^{-m}}( \cdo , y)|  \|^2_{L^2(X)} \le C
\|J^{\alpha}_j \|^2_{L^\infty}
 \mu(B(y,2^j r))^{-1}, 
\end{eqnarray*}
where
$
J_j^{\alpha}(r\lambda)=
(1+2^{2j}r^2\lambda^2)^m (r\lambda)^{2\alpha}F_j^{\alpha}(r\lambda)
$. Now to prove that the sum in (\ref{summ}) is bounded it is 
enough to show the following elementary estimate
\begin{equation}\label{pom}
\|J_j^{\alpha}(\lambda)\|_{L^\infty}=
\sup_{\lambda\in \R}|(1+2^{2j}\lambda^2)^m
(\lambda)^{\alpha}F_j^{\alpha}(\lambda)| \le 2^{-j}.
\end{equation}
As we noted before for any  fixed natural number $K\in \N$ we may assume that
$\Phi^{(l)}(0)=0$ when $1 \ge l \ge K$. Now  for any natural number $N$ 
we choose $K$ large enough
so that the Fourier transform of $H^{\alpha} (1 - \Phi_{1})$ 
is in   $C^\infty(\R-\{0\})$ and it
has a polynomial
decay of order $N$. Next  
$\widehat{F_j^{\alpha}}(\lambda)=
{[H^{\alpha} (1 - \Phi_{1})]}\widehat{{\,\,\,}}(\lambda)
(1-\varphi(2^{-j}\lambda))$ and it is not difficult to note that for 
any nonnegative integers $m_1$ and $m_2$ we can choose $K$ large enough
so that
$
\sup_{\lambda\in \R}|\lambda^{2m_1}F_j^{\alpha}(\lambda)| \le 2^{-m_2j}$
and (\ref{pom}) follows.

\end{proof}

\section{Hodge Laplacian and Riesz transform for $p>2$}\label{rhl}

Let $L_k$ be the Hodge-Laplace operator acting on 
$ L^2(\Lambda^kT^*M)$, where $M$ is $n$-dimension complete Riemannian 
manifolds. That is 
$$
\langle L_k \omega, \omega \rangle =
\int_M \big[ |d_k\omega|^2 
+|d_{k-1}^*\omega|^2\big]\dd \mu(x),
$$ 
where $d_k^*$ is the adjoint operator of $d_k$ (when $k=-1$ or $k=n$ one
should interpret $d_k\omega$ and $d_k^*\omega$ as $0$).
Note that if $\dd \mu=\nu(x) \dd x$, where $\dd x $ is the Riemannian measure
then 
\begin{equation}\label{dst}
d^*\omega \dd \mu= (-1)^{n(k+1)+1}*d*(\nu\omega)\dd x  
=(-1)^{n(k+1)+1}*(d*\omega+\frac{d\nu}{\nu} \wedge *\omega)\dd\mu
\end{equation} 
(see
\cite[(2) p. 220]{Wo}). 
Hence $L_k$
is defined by (\ref{hod1}) with $\beta_*=\frac{d\nu}{\nu}$.
Thus as a straightforward consequence of Theorems~\ref{fsp}~and~\ref{gl2} we
obtain the following Corollary
\begin{coro}\label{co}
Let $L_k$ be  the Hodge-Laplace operator acting on 
$ L^2(\Lambda^kT^*M)$, where $M$ is $n$-dimension complete Riemannian 
manifolds.
Suppose that for some number $N\in \N$ and points $x,y \in X$
there exist functions $V_x,V_y \colon \R^+ \mapsto \R$ such that
\begin{equation*} 
\| |K_{(I+t^2L_k)^{-N/4}}(z,\cdo)|  \|_{L^2(M)} \le  V_{z}(t)
    \quad \forall t > 0, z=x,y.
    \end{equation*}
Then, there exists a constant $C$ such that for all $t<\rho(x,y)^2$
\begin{eqnarray*}
|K_{\exp(-tL_k)}(x,y)| \le C  
V_x\bigg(\frac{t}{\d(x,y)}\bigg)V_y\bigg(\frac{t}{\d(x,y)}\bigg)
 \frac{{\exp}\left(\frac{-\d(x,y)^2}{4t}\right)}{\d(x,y)t^{-1/2}}.
\end{eqnarray*}
\end{coro}
{\em Remark.} Note that our proof of Theorem~\ref{gl}, contrary to some other
available arguments (see \cite{Co, DP, Gr1}), does not use the positivity of 
the heat kernel. Hence it works
for operator acting on fiber bundles and for the Hodge Laplacian in
particular.

As we mentioned in introduction one  cannot expect in general  
the Riesz transform
$d L_0^{-1/2}$ to be bounded on $L^p$ for $p>2$. However, the following 
theorem
shows that if $K_{(\exp(-tL_1)}(x,y)$ satisfies the expected on-diagonal bounds
then the Riesz transform
$d L_0^{-1/2}$ is bounded for all $p\in [2,\infty)$.
\begin{theorem}\label{pw2}
Suppose that $M$ is a complete Riemannian manifold satisfying the doubling
 volume
property i.e. Assumption~\ref{1a} and that $L_1$ satisfies
condition~(\ref{ell1}).
Then the Riesz transform $d_0 L_0^{-1/2}$ is bounded from $L^p(M)$ 
to $L^p(\Lambda^1T^*M)$ for
all $ 2 \le p<\infty$. 
\end{theorem} 
\begin{proof}
First we note that $d_0L_0=d_0d_0^*d_0=d^*_1d_1d_0+d_0d_0^*d_0=L_1d_0$ so
 $d_0^*L_1 = L_0d_0^*$ and 
\begin{equation}
d_0^*{L_1}^{-1/2} = L_0^{-1/2}d_0^*.\label{com}
\end{equation}
Hence $(d_0{L_0}^{-1/2})^*= L_0^{-1/2}d_0^*=d_0^*{L_1}^{-1/2}\!.$
 Now $d_0L_0^{-1/2}$ is bounded from $L^p(M)$  to $L^p(\Lambda^1T^*M)$ 
for some $2 \le p < \infty$ if 
and only if $(d_0{L}^{-1/2})^*= d_0^*{L_1}^{-1/2}$
 is bounded  
from $L^{p'}(\Lambda^1T^*M)$  to $L^{p'}(M)$ for $1/p+1/p'=1$. 
However, by (\ref{dst}) the operator $d_0^*$ is local and we assume
 that ${L_1}$ satisfies
(\ref{ell1}) so continuity of 
$d_0^* {L_1}^{-1/2}$ on $L^{p'}$ for $1<p' \le 2$ 
follows from Theorems~\ref{gl}~and~\ref{fsp}.
\end{proof}
Theorem~\ref{pw2} generalizes results described in \cite{CD2, CD3} to a very 
natural and somehow optimal setting.

\section{Other applications}

\subsection{Schr\"odinger operators}\label{sh}

Let $M$ be a connected and complete Riemannian manifold.
The Riemannian metric give us canonical isomorphisms 
$\Lambda^1T_x^*M \cong \Lambda^1T_xM$. We denote this isomorphism by
\emph{tilde}, so if $\omega \in \Lambda^1T_x^*M$ then $\widetilde{\omega}$ is
the corresponding dual element in $\Lambda^1T_xM$
and if $Y \in \Lambda^1T_xM$ then $\widetilde{Y}$ is
the corresponding dual element in $\Lambda^1T_x^*M$ . Then if $f$ is a function
on $M$, its gradient is the vector $\gr f=\widetilde{df}$.
We consider the operator $L_{Y,V}$ given by the formula
\begin{eqnarray}\label{ujaa}
\langle L_{Y,V}f,f\rangle =\int\limits_{M}\,\big(|\gr f(x)+if(x)Y|^2%
+V^2(x)|f(x)|^2\big)\dd \mu(x)\\=
\int\limits_{M}\,\big(|df+if\widetilde{Y}|^2
+|*f(x)\wedge V |^2\big)\dd \mu(x),
\end{eqnarray}
where $f \in C^{\infty}_c(M)$,
$Y$ is   a real   vector  field  such  that  $|Y|^2 \in
L^1_{\mbox{\small{loc}}}(M)$,     $V  \in
L^2_{\mbox{\small{loc}}}(M)$.   
\begin{theorem}\label{cde}
Suppose that the manifold $M$ satisfies Assumption~\ref{1a} 
and that the operator $L_{Y,V}$
satisfies (\ref{ell2}).
Then the operators
\begin{equation*}
VL_{Y,V}^{-1/2} \quad \mbox{and}  \quad ( \gr -iY)L_{Y,V}^{-1/2}
\end{equation*}
are bounded on $L^p$ for all $p\in (1,2]$ and 
of weak type $(1,1)$. 
\end{theorem}
\begin{proof}
Theorem~\ref{cde} is a straightforward consequence of 
Theorems~\ref{gl}~and~\ref{fsp}.
\end{proof}
{\em Remarks.} 1. In the case $M=\R^n$ and $L_{0,0}=\Delta$ 
Theorem~\ref{cde} was obtained
independently in \cite{CDE}.  In the same setting, i.e  $M=\R^n$ and 
$L_{0,0}=\Delta$ 
the Riesz transform was studied in
\cite{Sh1, Sh2}. Theorem~\ref{cde} applied to this setting yields an
interesting variation of Shen's result \cite[Theorem~0.5]{Sh1}
without any assumption concerning regularity of the potential $V$.
The counterexample investigated in \cite{Sh1},
($V(x)=|x|^{2-\varepsilon}$) shows that the operator
$\nabla (\Delta+V)^{-1/2}$ is not necessarily bounded for $p>2$. 
However, boundedness of the Riesz transforms for a larger range of $p$ can
be obtained if one imposes an additional regularity conditions for the potential
$V$ (see again \cite{Sh1}). 

2. It follows from \cite[Theorem~2.3]{Sim} and \cite[Theorem~{4.2}, p. 470]{1}
that 
\begin{equation}\label{Ar}
|K_{\exp{-tL_{Y,V}}}(x,y)| \le K_{\exp{-tL_{0,0}}}(x,y).
\end{equation}
Hence if $L_{0,0}$  satisfies
\ref{ell2} (or \ref{ell1}) then $L_{Y,V}$ also  satisfies this
assumption.

\subsection{Sub-elliptic operators acting on Lie groups}
Now let me describe  another application of Theorem~\ref{gl}.
Let $G$ be a Lie group with polynomial growth. For a system of
left-invariant vector fields $X_1,\ldots,X_k$ satisfying H\"ormander
condition, a function  $V  \in
L^2_{\mbox{\small{loc}}}(G)$ and a family of
functions $Y_1,\ldots,Y_k \in L^2(G)$ we define the operator $L_{V,Y}$ by the
formula
$$
L_{V,Y}=-\sum_{j=1}^k(X_j-iY_j)^2 +|V|^2.
$$
One can easily notice that the proof of  Theorem~\ref{fsp} works for
the operator $L_{V,Y}$. Hence $L_{V,Y}$ satisfies the finite speed 
propagation theorem with respect to the optimal control metrics corresponding
to the system $X_1,\ldots,X_k$ (see e.g.  \cite[\S III.4, p. 39]{VSC} for the definition).  
And so the  following theorem is again a straightforward consequence of 
Theorem~\ref{gl}.
\begin{theorem}\label{sie}
The Riesz transforms
\begin{equation}\label{lii}
VL_{V,Y}^{-1/2} \quad \mbox{and}  \quad ( X_j-iY_j)L_{V,Y}^{-1/2} \quad
\mbox{for}
\quad j=1,\ldots, k
\end{equation}
are bounded on $L^p$ for all $p\in (1,2]$ and 
of weak type $(1,1)$. 
\end{theorem}
Theorem~\ref{sie} is related to results described in \cite{Li}.
In \cite[Theorem C]{Li} Li 
proved the continuity  of the Riesz transforms (\ref{lii})
with additional assumptions about the regularity of the potential
$V$ in the case $Y_j=0$ for all $j$.

\end{document}